\documentclass[11pt,reqno,a4paper,oneside]{amsart}

\usepackage{latexsym}
\usepackage{amsmath}
\usepackage{amsthm}
\usepackage{amssymb}
\usepackage{amsfonts}
\usepackage{fancyhdr}
\usepackage{comment}

\usepackage{mathrsfs}                           % For different calligraphic fonts

\newcommand{\mR}{\mathbf{R}}                    % Formatting for R
\newcommand{\abs}[1]{\lvert #1 \rvert}          % Formatting for the absolute value
\newcommand{\norm}[1]{\lVert #1 \rVert}         % Formatting for the norm
\newcommand{\br}[1]{\langle #1 \rangle}         % Formatting for the inner product
\newcommand{\mS}{\mathscr{S}}

\newcommand{\mSp}{\mathscr{S}^{\prime}}

\newcommand{\ehat}{\,\hat{\rule{0pt}{6pt}}\,}
\newcommand{\edot}{\,\dot{\rule{0pt}{6pt}}\,}

\newcommand{\supp}{\mathrm{supp}}

\theoremstyle{definition}
\newtheorem{thm}{Theorem}[section]

\newtheorem{lemma}{Lemma}[section]

\newtheorem*{remark}{Remark}

\title[Stability for wave equations]{Stability for Solutions of Wave Equations \\
with $C^{1,1}$ Coefficients}

\author{Mikko Salo}
\address{Department of Mathematics and Statistics / RNI,
University of Helsinki}
\email{mikko.salo\@@helsinki.fi}

\begin{document}

\begin{abstract}
We consider the stable dependence of solutions to wave equations on metrics in $C^{1,1}$ class. The main result states that solutions depend uniformly continuously on the metric, when the Cauchy data is given in a range of Sobolev spaces. The proof is constructive and uses the wave packet approach to hyperbolic equations.
\end{abstract}

\maketitle

\section{Introduction}

We consider the wave equation in $\mR_t \times \mR^n_x$,
\begin{equation*}
\left\{ \begin{array}{rl}
(D_t^2 - A(x,D_x))u(t,x) &\!\!\!= F(t,x), \\[4pt]
u|_{t=0} &\!\!\!= f, \\[4pt]
\partial_t u|_{t=0} &\!\!\!= g.
\end{array} \right.
\end{equation*}
Here $A(x,D_x) = a^{ij}(x) D_{x_i} D_{x_j}$ is a uniformly elliptic operator, satisfying $a^{ij} = a^{ji}$, and $a^{ij}(x) \xi_i \xi_j \geq c \abs{\xi}^2$ for $\xi \in \mR^n$. We assume that the functions $a^{ij}$ are in the space $C^{1,1}(\mR^n)$, with norm $\norm{a}_{C^{k-1,1}} = \sum_{\abs{\alpha} \leq k} \norm{\partial^{\alpha} a}_{L^{\infty}}$.

The question investigated in this article is the stable dependence of the solution $u$ on the metric $(a^{ij})$. Intuition for the problem can be obtained from the simplest possible case, namely the one-dimensional wave equation with constant sound speed $c$. Given $f \in L^2(\mR)$, the equation 
\begin{eqnarray*}
 & (\partial_t^2 - c^2 \partial_x^2) u(t,x) = 0, & \\
 & u(0) = f, \partial_t u(0) = 0 &
\end{eqnarray*}
has the solution $u(t,x) = \frac{1}{2} [ f(x-ct) + f(x+ct) ]$. Thus, the solution at time $t = 1$ is obtained by translating $f$ by $c$ units in the positive and negative directions. Since translation of $L^2$ functions is a uniformly continuous operation, we see that $u(1,\,\cdot\,)$ depends uniformly continuously in $L^2$ norm on the sound speed $c$. A stronger result may be obtained if the initial data is smoother: if $f \in H^1$, then $\norm{f(\,\cdot\,-c) - f(\,\cdot\,-c')}_{L^2} \leq \norm{\nabla f}_{L^2} \abs{c-c'}$, and $u(1,\,\cdot\,)$ depends Lipschitz continuously in $L^2$ on the sound speed.

Our main result is the following theorem, which shows that uniformly continuous or Lipschitz dependence are valid also for general wave equations with $C^{1,1}$ metrics. In this introduction, we state the theorem only in the case where an initial velocity $g$ is present. In Section 7 we will give the straightforward extensions to cases where driving terms $F$ and initial positions $f$ are present, and also where $A(x,D_x)$ is replaced by a divergence form or Laplace-Beltrami operator.

\begin{thm} \label{thm:main}
Suppose $M$ is a large constant such that 
\begin{eqnarray}
 & \norm{a^{ij}}_{C^{1,1}} \leq M, \quad a^{ij} \xi_i \xi_j \geq M^{-2} \abs{\xi}^2, & \label{a_assumptions} \\
 & t \in [-M,M]. & \label{t_assumptions}
\end{eqnarray}
If $-1 \leq \alpha \leq 2$, then for each $g \in H^{\alpha}$ there is a unique weak solution $u \in C([-M,M];H^{\alpha+1}) \cap C^1([-M,M];H^{\alpha})$ of the problem 
\begin{equation} \label{main_cauchy}
\left\{ \begin{array}{rl}
(D_t^2 - A(x,D_x))u(t,x) &\!\!\!= 0, \\[4pt]
u|_{t=0} &\!\!\!= 0, \\[4pt]
\partial_t u|_{t=0} &\!\!\!= g.
\end{array} \right.
\end{equation}
If $-1 \leq \alpha < 2$, and if $A = (a^{ij})$, $B = (b^{ij})$ satisfy \eqref{a_assumptions} and $u_A$, $u_B$ are the corresponding solutions, then for any $\varepsilon > 0$ there is $\delta > 0$ such that 
\begin{equation} \label{uniform_stability}
\norm{u_A(t) - u_B(t)}_{H^{\alpha+1}} < \varepsilon \ \text{ whenever }\ \norm{A-B}_{C^{0,1}} < \delta.
\end{equation}
Further, if $-1 \leq \alpha \leq 1$ and $g \in H^{\alpha+1}$, then 
\begin{equation} \label{lipschitz_stability}
\norm{u_A(t) - u_B(t)}_{H^{\alpha+1}} \leq C \norm{A-B}_{C^{0,1}} \norm{g}_{H^{\alpha+1}}
\end{equation}
where $C$ depends only on $M$ and $n$.
\end{thm}

Here, by a weak solution we mean a function $u \in C([-M,M];H^{\alpha+1}) \cap C^1([-M,M];H^{\alpha})$ which solves the equation in the sense of distributions, and satisfies the initial conditions in the vector-valued sense.

Existence and uniqueness of solutions to linear wave equations is of course classical and can be established via energy estimates under quite general conditions. We refer to \cite{stolk} for a comprehensive account and further references. A constructive method, valid in the setting of our theorem, for proving existence and uniqueness of weak solutions was introduced in \cite{smith1} using a wave packet approach.

Stability estimates such as the ones in Theorem \ref{thm:main} are also classical, see \cite{stolk} and references given there. In Section 7 we give an easy argument which uses just the existence and uniqueness of solutions, without any precise knowledge about the solutions.

The novelty is that our proof of Theorem \ref{thm:main} is constructive: an explicit expression for the solution is given, and the stability properties are deduced from that. The constructive method gives the same intuition to stability as in the $n=1$ case, namely that stability for solutions should be the same as for translating functions. We use the wave packet approach introduced in \cite{smith1}, where the initial data is decomposed into wave packets, and an approximate solution to the equation is obtained by translating wave packets along the Hamilton flow. The stability is exactly governed by this translation.

The main motivation for this study comes from inverse problems in seismic imaging, where the analysis of solutions of wave equations has a key role. Many existing results, see \cite{dehoop_msri} for a survey, assume a linearization about a smooth sound speed $c_0$ and use microlocal analysis and calculus of Fourier integral operators (FIOs).

There has been recent interest, see \cite{dsu}, in the practically more realistic case where $c_0$ is not smooth, and in this case few results are known. One reason is that there is no calculus of nonsmooth FIOs. However, solution operators for nonsmooth wave equations are understood quite well due to the wave packet approach of \cite{smith1}, and a more precise analysis of these operators is expected to lead to new results. This work is an attempt in this direction. The related work \cite{dsu} discusses propagation of singularities (in terms of concentration of wave packets) and application of wave packets in numerical computations. We remark that wave packets are the same as curvelets \cite{candesdemanet}, \cite{candesdonoho}, which have been introduced in image processing as an efficient way of representing functions with singularities on smooth curves.

The proof of Theorem \ref{thm:main} is based on constructing an explicit solution operator for the nonsmooth wave equation, following the method of \cite{smith1}. The idea is to localize the initial data to dyadic frequency shells, and write each localized piece as a superposition of wave packets at the given frequency. The action of the wave group on a wave packet is well approximated by translating the packet along the Hamilton flow, and this gives an approximate solution operator for frequency localized initial data. These are added up to obtain a parametrix for the full equation.

Actually, to handle the nonsmooth coefficients, at frequency level $2^k$ one truncates the coefficients to frequencies less than $2^{k/2}$ and uses the Hamilton flow for the truncated metric. The error terms resulting from this will add up to a bounded operator, and this error can be iterated away by solving a Volterra equation.

The precise construction of solution operator will be a combination of methods in \cite{smith1} and \cite{smith2} and it does not involve new ideas. However, in view of the stability result, we need to give the construction in detail to see how the operator depends on the metric. The main outline is the same as in \cite{smith1}, which used a discrete wave packet frame. Some computations are simplified if one uses instead a continuous wave packet representation (i.e.~the FBI transform, see \cite{delort}). This was used in \cite{tataru1}, \cite{tataru2}, \cite{tataru3} for wave packets based on the Gaussian. We will follow \cite{smith2} which used instead wave packets compactly supported in frequency, a property which will keep the dyadic frequency annuli separated.

The crux of the stability proof is Lemma \ref{flowstability}, which considers the stability of translating along Hamilton flow. Lipschitz stability for translation involves a loss of one derivative, and the main point in the proof is that there is a smooth deformation of two Hamilton flows obtained by deforming the corresponding metrics. For the full result, one needs to check that the several corrections required to obtain an exact solution operator, most importantly the Volterra iteration, do not affect the stability given by translation.

We first prove Lipschitz stability with a loss of one derivative, and the uniform continuity is an immediate consequence. Since there is a loss of one derivative arising from the flow, we can afford to lose one derivative in certain other estimates as well. In this way, we get stability in terms of the $C^{0,1}$ norm of the metric instead of $C^{1,1}$ norm.

The plan of the paper is as follows. Section 2 contains some basic facts about Hamilton flows, and Section 3 introduces the FBI transform. In Section 4 we outline the construction of the solution operator, and Section 5 contains the details. The stability result, Theorem \ref{thm:main}, is proved in Section 6. Section 7 discusses variations of the stability result.

\vspace{5pt}

\noindent {\bf Notation.} We write $D_t = \frac{1}{i} \partial_t$ and $D_{x_j} = \frac{1}{i} \partial_{x_j}$. The gradient with respect to $x$ or $\xi$ is denoted by $d_x$ or $d_{\xi}$, and $D = \frac{1}{i} d$. Throughout the paper, $M$ will be a large constant so that \eqref{a_assumptions} and \eqref{t_assumptions} are satisfied. We write $a \lesssim b$ if $a \leq Cb$ where $C > 0$ depends only on $M$ and the dimension $n$. Also, $a \sim b$ means that $a \lesssim b$ and $b \lesssim a$. We write $L^p = L^p(\mR^n)$, similarly for the $L^2$ Sobolev spaces $H^{\alpha}$, the spaces $C^{0,1}$ and $C^{1,1}$, and the Schwartz space $\mS$. The mixed norm spaces are denoted by $L^p_t H^{\alpha}_x = L^p([-M,M] ; H^{\alpha})$, similarly $C^j_t H^{\alpha}_x = C^j([-M,M] ; H^{\alpha})$.

\subsection*{Acknowledgements}

Research partly supported by the Academy of Finland. Part of this research was carried out during visits at MSRI and at the University of Washington, and I wish to express my gratitude to these institutions for their hospitality and support. I would like to thank Hart Smith for generous advice, and Gunther Uhlmann for suggesting the problem.

\section{Hamilton flow}

We record for later use some elementary facts related to Hamilton flows. Let $\chi(\xi)$ be a smooth cutoff supported in the unit ball with $\chi = 1$ for $\abs{\xi} \leq 1/2$. We define smooth approximations of the coefficients $a^{ij}$ by 
\begin{equation*}
a^{ij}_k(x) = \chi(2^{-k/2}D_x) a^{ij}(x).
\end{equation*}
Then $a^{ij}_k$ is supported in frequency in $\{\abs{\xi} \leq 2^{k/2}\}$, and satisfies derivative bounds 
\begin{equation} \label{truncation_est}
\norm{\partial^{\alpha} a^{ij}_k}_{L^{\infty}} \lesssim 2^{\frac{k}{2} \max(0,\abs{\alpha}-2)}. \\[5pt]
\end{equation}

Consider the Hamilton flow related to $D_t^2 - A_k(x,D_x)$, where $A_k(x,D_x) = a^{ij}_k(x) D_{x_i} D_{x_j}$. It will be useful to define this in terms of the functions 
\begin{equation*}
p_k^{\pm}(x,\xi) = \pm \chi(2^{-k/2} D_x)\sqrt{A_k(x,\xi)}.
\end{equation*}
The Hamilton equations are given by 
\begin{align*}
\dot{x}(t) &= d_{\xi} p(x(t),\xi(t)), \\
\dot{\xi}(t) &= -d_x p(x(t),\xi(t)),
\end{align*}
where $p = p_k^{\pm}$, $x = x_k^{\pm}$, and $\xi = \xi_k^{\pm}$. Here $(x,\xi) \mapsto (d_{\xi} p(x,\xi), -d_x p(x,\xi))$ is a smooth vector field in $T^* \mR^n = \mR^n_x \times (\mR^n_{\xi} \smallsetminus \{0\})$, so given an initial condition $(x(0),\xi(0)) = (y,\eta) \in T^* \mR^n$, the Hamilton equations have a solution $(x(t),\xi(t))$ depending smoothly on $t,y,\eta$ at least for small time.

It is well known that the solution $(x(t),\xi(t))$ exists for all time. For, if it exists for $t$ in an interval $I = (-t_0,t_0)$, then $\xi(t) \neq 0$ for $t \in I$, and $h(t) = \abs{\xi(t)}^2$ satisfies by the homogeneity of $p$ 
\begin{equation*}
\dot{h}(t) = 2 \xi(t) \cdot \dot{\xi}(t) = -2 \xi(t) \cdot d_x p(x(t),\xi(t)) = m(t) h(t)
\end{equation*}
where $m(t) = -2 \omega(t) \cdot d_x p(x(t),\omega(t))$ and $\omega(t) = \xi(t)/\abs{\xi(t)}$. Solving the equation gives $h(t) = e^{\int_0^t m(s) \,ds} h(0)$. But $\abs{m(t)} \leq 2 \sup_{x \in \mR^n, \abs{\omega} = 1} \abs{d_x p(x,\omega)}$, and this argument shows that 
\begin{equation*}
e^{-t_0 \sup \abs{d_x p(x,\omega)}} \abs{\xi(0)} \leq \abs{\xi(t)} \leq e^{t_0 \sup \abs{d_x p(x,\omega)}} \abs{\xi(0)}
\end{equation*}
when $t \in I$. For $x(t)$ we have $\abs{\dot{x}(t)} \leq \sup_{x \in \mR^n, \abs{\omega} = 1} \abs{d_{\xi} p(x,\omega)}$ for $t \in I$, so $\abs{x(t) - x(0)} \leq t_0 n^{1/2} \sup \abs{d_{\xi} p(x,\omega)}$ for $t \in I$. Thus $(x(t),\xi(t))$ stays in a compact subset of $T^* \mR^n$ when $t \in I$, and we may extend the solution past the endpoints of $I$ and to all time.

The following consequence will be used frequently: if $(x(t),\xi(t))$ satisfy the Hamilton equations and $\abs{\xi(0)} \sim \lambda$, then 
\begin{equation*}
\abs{\xi(t)} \sim \lambda
\end{equation*}
for $\abs{t} \leq M$.

We write $\chi_{t,s} = \chi_{t,s}^{k,\pm}$ for the canonical transformation $(y,\eta) \mapsto (x(t),\xi(t))$, where $(x,\xi)$ solve the Hamilton equations with initial condition $(x(s),\xi(s)) = (y,\eta)$. Then $\chi_{t,s}$ is a symplectic diffeomorphism of $T^* \mR^n$ with inverse $\chi_{s,t}$.

\section{Wave packet representation}

Let $g$ be a real, even Schwarz function in $\mR^n$ with $\norm{g}_{L^2} = (2\pi)^{-\frac{n}{2}}$, and assume $\hat{g}$ is supported in the unit ball. For $\lambda \geq 1$ and $y,x,\xi \in \mR^n$ define 
\begin{equation} \label{glambda_def}
g_{\lambda}(y ; x,\xi) = \lambda^{\frac{n}{4}} e^{i\xi \cdot (y-x)} g(\lambda^{\frac{1}{2}}(y-x)).
\end{equation}
This is a wave packet at frequency level $\lambda$, centered in space at $x$ and in frequency at $\xi$. The Fourier transform is given by 
\begin{equation} \label{glambdahat_def}
\hat{g}_{\lambda}(\eta ; x,\xi) = \lambda^{-\frac{n}{4}} e^{-i\eta \cdot x} \hat{g}(\lambda^{-\frac{1}{2}}(\eta-\xi)).
\end{equation}
The FBI transform of a function $f \in \mS(\mR^n)$ is given by 
\begin{equation*}
T_{\lambda} f(x,\xi) = \int f(y) \overline{g_{\lambda}(y;x,\xi)} \,dy.
\end{equation*}
Suppose $\lambda \geq 2^6$. Then if $\hat{f}$ is supported in $\frac{1}{4} \lambda < \abs{\xi} < \lambda$, then $T_{\lambda} f$ vanishes unless $\frac{1}{8} \lambda < \abs{\xi} < 2 \lambda$. If $F \in \mS(\mR^{2n}_{x,\xi})$, the adjoint has the form 
\begin{equation*}
T_{\lambda}^* F(y) = \iint F(x,\xi) g_{\lambda}(y ; x,\xi) \,dx\,d\xi.
\end{equation*}
It follows that $T_{\lambda}^* T_{\lambda} = I$, and $\norm{T_{\lambda} f}_{L^2(\mR^{2n}_{x,\xi})} = \norm{f}_{L^2(\mR^n)}$.

The following result, stating the $L^2$ boundedness of FBI transform type operators, is from \cite[Lemma 3.1]{smith2}.

\begin{lemma} \label{fbisobolevbounded}
Let $g_{x,\xi}$ be a family of Schwartz functions whose Schwartz seminorms are bounded uniformly in $x$ and $\xi$. Then the operator $T$, defined for Schwartz functions by 
\begin{equation*}
Tf(x,\xi) = \int f(y) \overline{(g_{x,\xi})_{\lambda}(y; x,\xi)} \,dy
\end{equation*}
is bounded from $L^2(\mR^n)$ to $L^2(\mR^{2n})$, and the adjoint $T^*$ given by 
\begin{equation*}
T^* F(y) = \iint F(x,\xi) (g_{x,\xi})_{\lambda}(y;x,\xi) \,dx\,d\xi
\end{equation*}
is bounded from $L^2(\mR^{2n})$ to $L^2(\mR^n)$. The norms of $T$ and $T^*$ are bounded by a sum of finitely many Schwartz seminorms of the $g_{x,\xi}$.
\end{lemma}

\section{Outline of construction of solution operator} \label{sec:parametrixoutline}

We will outline the construction of a solution operator $S(t): g \mapsto u(t,\,\cdot\,)$ for the problem \eqref{main_cauchy}. More details are given in Section \ref{sec:parametrixdetails}.

Let us start with the standard Littlewood-Paley frequency localization. Let $\chi(\xi)$ be a smooth cutoff supported in the unit ball with $\chi = 1$ for $\abs{\xi} \leq 1/2$, and take $\beta_k(D)$ to be a Littlewood-Paley partition of unity with 
\begin{equation*}
\beta_0(\xi) + \sum_{k=1}^{\infty} \beta_k(\xi) = 1
\end{equation*}
where $\beta_0$ is supported in the unit ball, $\beta_1$ is supported in $\{ \frac{1}{2} < \abs{\xi} < 2 \}$, and $\beta_k(\xi) = \beta_1(2^{-k+1} \xi)$.

We let $g_k = \beta_k(D)g$, and for each $k$ we want to find an approximate solution of 
\begin{equation} \label{localizedwaveeq}
\left\{ \begin{array}{rl}
(D_t^2 - A_k(x,D_x))u_k &\!\!\!= 0, \\[4pt]
u_k(0) = 0,\ \partial_t u_k(0) &\!\!\!= g_k.
\end{array} \right.
\end{equation}
This will be done by reducing matters to first order hyperbolic equations. Consider the pseudodifferential symbol 
\begin{equation*}
p_k^{\pm}(x,\xi) = \pm \chi(2^{-k/2}D_x) \sqrt{A_k(x,\xi)}
\end{equation*}
so that $p_k^{\pm}$ is positively homogeneous of degree $1$ in $\xi$ and smooth when $\xi \neq 0$, and satisfies for $\abs{\xi} = 1$
\begin{equation*}
\abs{\partial_x^{\alpha} \partial_{\xi}^{\beta} p_k^{\pm}(x,\xi)} \lesssim 2^{\frac{k}{2} \max(0,\abs{\alpha}-2)}.
\end{equation*}
Thus $p_k^{\pm} \beta_k(\xi) \in S^1_{1,1/2}$. The Fourier transform in $x$, $(p_k^{\pm})\ehat(\,\cdot\,,\xi)$, is supported in $\{ \abs{\eta} \leq 2^{k/2}\}$. Since 
\begin{equation*}
(p(x,D) f) \ehat(\eta) = (2\pi)^{-n} \int \hat{p}(\eta-\xi,\xi) \hat{f}(\xi) \,d\xi,
\end{equation*}
we see that $p_k^{\pm}(x,D) f$ will be supported in $\abs{\eta} \sim 2^k$ if $\hat{f}$ is supported in $\abs{\eta} \sim 2^k$. We will also use the symbols 
\begin{equation*}
q_k^{\pm} = \chi(2^{-k/2}D_x) (1/p_k^{\pm})
\end{equation*}
which are homogeneous of degree $-1$ in $\xi$, and $q_k^{\pm} \beta_k(\xi) \in S^{-1}_{1,1/2}$.

The approximate solution operators related to \eqref{localizedwaveeq} are 
\begin{equation*}
E_k^{\pm}(t) g = T_k^* U_k^{\pm}(t) T_k( \frac{i}{2} Q_k^{\pm} \beta_k(D)g)
\end{equation*}
where $T_k = T_{\lambda}$ with $\lambda = 2^k$, and $U_k^{\pm}(t)$ is translation for time $-t$ along the Hamilton flow of $p_k^{\pm}$, that is, 
\begin{equation*}
U_k^{\pm}(t)F = F \circ \chi_{0,t}^{k,\pm}.
\end{equation*}
We are ready to write down the first approximate solution operator for \eqref{main_cauchy}. It will be given by 
\begin{equation*}
\tilde{S}(t) g = t \sum_{k<k_0} g_k + \sum_{k \geq k_0} (u_k^+ + u_k^-)
\end{equation*}
where $k_0$ is a constant depending on $M$ which will be chosen later, and where 
\begin{equation*}
u_k^{\pm} = E_k^{\pm}(t) g.
\end{equation*}

Since $Q_k^{\pm}$ is of order $-1$, it is not hard to see that $\tilde{S}(t)$ is an operator of order $-1$, i.e.~it maps $H^{\alpha}(\mR^n)$ to $H^{\alpha+1}(\mR^n)$. The main point is that the operator 
\begin{equation*}
\tilde{R}_k^{\pm}(t) = (D_t + P_k^{\pm}(y,D_y)) E_k^{\pm}(t),
\end{equation*}
which is a half-wave operator applied to $E_k^{\pm}(t)$, is of order $-1$. This implies that when one applies the wave operator $D_t^2 - A(x,D_x)$ to $\tilde{S}(t)$, the resulting operator will be of order $0$, which is one order lower than expected. Therefore we really have an approximate solution operator.

However, since $\partial_t \tilde{S}(t) g|_{t = 0}$ is only close to $g$ but not equal to $g$ in general, we need to correct the initial values of the operator. Thus instead of $\tilde{S}(t)$ we use 
\begin{equation*}
\widehat{S}(t) = \tilde{S}(t)(I+K)^{-1}
\end{equation*}
where $K$ is given by 
\begin{equation*}
K = \sum_{k \geq k_0} (i \tilde{R}_k^+(0) + i \tilde{R}_k^-(0) + \frac{1}{2} R_k^+ \beta_k(D) + \frac{1}{2} R_k^- \beta_k(D))
\end{equation*}
Here $R_k^{\pm} \beta_k(D)$ are order $-1$ operators given by 
\begin{equation*}
P_k^{\pm} Q_k^{\pm} \beta_k(D) = (I + R_k^{\pm}) \beta_k(D).
\end{equation*}
If $k_0$ is chosen large enough then $K$ will be small on $H^{\alpha}(\mR^n)$, and $\widehat{S}(t)$ will be an operator of order $-1$ with $\widehat{S}(t)g |_{t=0} = 0$ and $\partial_t \widehat{S}(t)g |_{t=0} = g$. Thus $\widehat{S}(t)g$ will be an approximate solution of \eqref{main_cauchy} with correct Cauchy data.

It remains to show that one may iterate away the error and obtain an exact solution. To do this we seek a solution of \eqref{main_cauchy} of the form 
\begin{equation*}
u(t,x) = \widehat{S}(t) g(x) + \int_0^t \widehat{S}(t,s) G(s,x) \,ds.
\end{equation*}
Here $\widehat{S}(t,s) = \widehat{S}(t-s)$ is the operator corresponding to $\widehat{S}(t)$ but where the initial surface is $\{t = s\}$ instead of $\{t = 0\}$. Then $\widehat{S}(t,s) g|_{t=s} = 0$ and $\partial_t \widehat{S}(t,s) g |_{t=s} = g$, and one has 
\begin{equation*}
\partial_t^2 \Big( \int_0^t \widehat{S}(t,s) G(s,x) \,ds \Big) = G(t,x) + \int_0^t \partial_t^2 \widehat{S}(t,s) G(s,x) \,ds.
\end{equation*}
We obtain 
\begin{equation*}
(D_t^2 - A(x,D_x)) u(t,x) = T(t)g(x) - G(t,x) + \int_0^t T(t,s) G(s,x) \,ds
\end{equation*}
where $T(t,s) = (D_t^2 - A(x,D_x)) \widehat{S}(t,s)$. As remarked above, we will show that $T(t,s)$ is an operator of order $0$, which in the present setting with a $C^{1,1}$ wave operator will mean that it is bounded on $H^{\alpha}(\mR^n)$ for $-1 \leq \alpha \leq 2$. Then $u$ will be a solution provided that $G(t,x) = V(T(t)g(x))$, where $G = VF$ solves the Volterra equation 
\begin{equation*}
G(t,x) - \int_0^t T(t,s) G(s,x) \,ds = F(t,x).
\end{equation*}
Since $T(t,s)$ is bounded on $H^{\alpha}$ also $V$ is bounded on $L^{\infty}_t H^{\alpha}_x$, for $-1 \leq \alpha \leq 2$. Thus the full solution operator for \eqref{main_cauchy} will be 
\begin{equation*}
S(t)g(x) = \widehat{S}(t)g(x) + \int_0^t \widehat{S}(t,s) V(T(s) g(x)) \,ds.
\end{equation*}
We will now give the details.

\section{Details of construction of solution operator} \label{sec:parametrixdetails}

Let $(a^{ij})$ be a symmetric matrix of $C^{1,1}$ functions satisfying \eqref{a_assumptions}. We take $M$ large enough so that \eqref{a_assumptions} holds also for the truncated metrics $(a^{ij}_k)$, and we also assume \eqref{t_assumptions}. We will use all the notations in Section \ref{sec:parametrixoutline}. Also, $k_0$ will be a sufficiently large integer, depending on $M$ and $n$.

We start by noting that if the Cauchy data is localized near frequency $2^k$, then the operators constructed in the preceding section preserve this localization.

\begin{lemma}
Suppose the Fourier transform of $f$ is supported in $\abs{\xi} \sim 2^k$. Then $T_k f$, $P_k^{\pm} f$, $Q_k^{\pm} f$ vanish unless $\abs{\xi} \sim 2^k$, provided $k \geq k_0$. If $F(x,\xi)$ vanishes unless $\abs{\xi} \sim 2^k$, then $T_k^* F$ and $U_k^{\pm}(t) F$ vanish unless $\abs{\xi} \sim 2^k$, provided $k \geq k_0$.
\end{lemma}

Using this lemma, it will be enough to consider a fixed frequency and the result will follow by summing over dyadic annuli. Thus, let $\lambda = 2^k$, and write $p = p_k^{\pm}$, $U(t) = U_k^{\pm}(t)$, etc. The main idea is that the wave evolution of a wave packet at frequency $\lambda$ can be well approximated by transport along the Hamilton flow. The transport operator will be $L_{x,\xi} = L_k^{\pm}$, given by 
\begin{equation*}
L_{x,\xi} = d_{\xi} p(x,\xi) \cdot D_x - d_x p(x,\xi) \cdot D_{\xi}.
\end{equation*}
It is easy to see that $U(t)F$ will solve the corresponding transport equation.

\begin{lemma}
If $F \in \mS(\mR^n_x \times \mR^n_{\xi})$ then $U(t)F$ satisfies 
\begin{equation*} %\label{transporteq}
\left\{ \begin{array}{rl}
(D_t + L_{x,\xi})U(t)F &\!\!\!= 0, \\[4pt]
U(t)F|_{t=0} &\!\!\!= F
\end{array} \right.
\end{equation*}
\end{lemma}
%\begin{proof}
%If $G(t,x,\xi) = U(t)F(x,\xi)$, then $G(t,\chi_{t,0}(x,\xi)) = F(x,\xi)$. Differentiating this with respect to $t$ gives 
%\begin{multline*}
%D_t G(t,\chi_{t,0}(x,\xi)) + d_{\xi} p(\chi_{t,0}(x,\xi)) \cdot D_x G(t,\chi_{t,0}(x,\xi)) \\
% - d_x p(\chi_{t,0}(x,\xi)) \cdot D_{\xi} G(t,\chi_{t,0}(x,\xi)) = 0.
%\end{multline*}
%Since $\chi_{t,0}$ is bijective on $T^* \mR^n$ and $\chi_{0,0} = \mathrm{Id}$, the result follows.
%\end{proof}

Let now $f$ be a function localized near frequency $\lambda$. We write $f$ as a superposition of wave packets at frequency $\lambda$, 
\begin{equation*}
f(y) = T_{\lambda}^* T_{\lambda} f(y) = \iint T_{\lambda} f(x,\xi) g_{\lambda}(y;x,\xi) \,dx\,d\xi.
\end{equation*}
To get an approximate solution $u$ to the half-wave equation $(D_t + P_y) u(t,y) = 0$ with $u(0,y) = f(y)$, where $P_y = p(y,D_y)$, we transport the wave packets for time $t$ along the Hamilton flow of $p$. Then $u$ is given by 
\begin{equation*}
u(t,y) = \iint T_{\lambda} f(x,\xi) g_{\lambda}(y;\chi_{t,0}(x,\xi)) \,dx\,d\xi.
\end{equation*}
Using that $\chi_{t,0}$ is a symplectic map with inverse $\chi_{0,t}$, $u$ will be equal to 
\begin{equation*}
u(t,y) = T_{\lambda}^* U(t) T_{\lambda} f.
\end{equation*}
To measure how far $u$ is from an exact solution, we compute 
\begin{multline*}
D_t u = T_{\lambda}^* D_t U(t) T_{\lambda} f = -T_{\lambda}^* L_{x,\xi} U(t) T_{\lambda} f \\
 = \iint (U(t) T_{\lambda} f)(x,\xi) L_{x,\xi} g_{\lambda}(y;x,\xi) \,dx\,d\xi,
\end{multline*}
the last equality by integration by parts. Thus 
\begin{equation} \label{halfwaveu}
(D_t + P_y) u = \iint (U(t) T_{\lambda} f)(x,\xi) (P_y + L_{x,\xi}) g_{\lambda}(y;x,\xi) \,dx\,d\xi.
\end{equation}
The following lemma, corresponding to \cite[Lemma 3.2]{smith2}, will be crucial.

\begin{lemma} \label{wavepacketparameter}
Suppose $\abs{\xi} \sim \lambda$. Then 
\begin{equation*}
(P_y + L_{x,\xi}) g_{\lambda}(y ; x,\xi) = (g_{x,\xi})_{\lambda}(y ; x,\xi)
\end{equation*}
where $g_{x,\xi}$ is a family of functions whose Fourier transform is supported in a ball of radius $2$ and the Schwartz seminorms of $g_{x,\xi}$ are uniformly bounded in $x,\xi$. In fact $g_{x,\xi}(z) = m_{x,\xi}(z,D_z) g(z)$ where 
\begin{equation} \label{mxxidef}
m_{x,\xi}(z,\zeta) = \int_0^1 (1-s) \partial_s^2 ( p(x + s \lambda^{-1/2} z, \xi + s \lambda^{1/2} \zeta) ) \,ds.
\end{equation}
Here $m_{x,\xi}$ satisfies symbol estimates uniform in $x$ and $\xi$, 
\begin{equation} \label{mxxiest}
\abs{\partial_z^{\alpha} \partial_{\zeta}^{\beta} m_{x,\xi}(z,\zeta)} \lesssim \br{z}^2, \qquad \abs{\zeta} \leq 2.
\end{equation}
\end{lemma}
\begin{proof}
We define $g_{x,\xi}$ by the relation $(g_{x,\xi})_{\lambda} = (P_y + L_{x,\xi}) g_{\lambda}$. Recalling the formula \eqref{glambda_def} for $g_{\lambda}$, it follows that 
\begin{multline*}
(g_{x,\xi})_{\lambda} = (P_y - id_{\xi} p \cdot d_x + i d_x p \cdot d_{\xi}) g_{\lambda} = (P_y + id_{\xi} p \cdot d_y - d_x p \cdot (y-x)) g_{\lambda} \\
 = (2\pi)^{-n} \int e^{iy \cdot \eta} (p(y,\eta) - d_{\xi} p(x,\xi) \cdot \eta - d_x p(x,\xi) \cdot (y-x)) (g_{\lambda})\ehat(\eta) \,d\eta \\
 = (2\pi)^{-n} \int e^{iy \cdot \eta} \Big[ p(y,\eta) - p(x,\xi) - d_{\xi} p(x,\xi) \cdot (\eta-\xi) - d_x p(x,\xi) \cdot (y-x) \Big] (g_{\lambda})\ehat(\eta) \,d\eta
\end{multline*}
since $\xi \cdot d_{\xi} p = p$ by homogeneity. By \eqref{glambdahat_def}, we have 
\begin{multline*}
g_{x,\xi}(z) = (2\pi)^{-n} \int e^{iz \cdot \zeta} \Big[ p(x+\lambda^{-1/2}z,\xi+\lambda^{1/2}\zeta) - p(x,\xi) - d_x p(x,\xi) \cdot \lambda^{-1/2} z \\
 - d_{\xi} p(x,\xi) \cdot \lambda^{1/2} \zeta \Big] \hat{g}(\zeta) \,d\zeta
\end{multline*}
which shows that $(g_{x,\xi}) \ehat = 0$ outside a ball of radius $2$. From Taylor's formula we see that the term in brackets is equal to $m_{x,\xi}(z,\zeta)$ given by \eqref{mxxidef}. If $(\tilde{x}, \tilde{\xi}) = (x+s\lambda^{-1/2}z, \xi+s\lambda^{1/2}\zeta)$, we compute 
\begin{equation} \label{mxxi_explicit}
\partial_s^2( p(\tilde{x}, \tilde{\xi}) ) = \sum_{j,k}[ \partial_{x_j x_k} p(\tilde{x}, \tilde{\xi}) \lambda^{-1} z_j z_k + \\
 \partial_{x_j \xi_k} p(\tilde{x}, \tilde{\xi}) z_j \zeta_k + \partial_{\xi_j \xi_k} p(\tilde{x}, \tilde{\xi}) \lambda \zeta_j \zeta_k ]
\end{equation}
which gives \eqref{mxxiest}.
\end{proof}

Now we can prove that the half-wave operator applied to the approximate solution operator gives an operator of order $0$.

\begin{lemma} \label{halfwaveapproximate}
Suppose $\hat{f}$ is supported in $\abs{\xi} \sim \lambda$. Then 
\begin{equation*}
\norm{(D_t + P_y) T_{\lambda}^* U(t) T_{\lambda} f}_{L^2} \lesssim \norm{f}_{L^2}.
\end{equation*}
\end{lemma}
\begin{proof}
If $\hat{f}$ is supported in $\abs{\xi} \sim \lambda$, the same holds for $U(t) T_{\lambda} f$. The result follows from \eqref{halfwaveu}, Lemma \ref{wavepacketparameter} and Lemma \ref{fbisobolevbounded}.
\end{proof}

Next we apply a second order wave operator $D_t^2 - P_y^2$ to the approximate solution operator, and show that there is a loss of one derivative.

\begin{lemma} \label{psdofullwaveapproximate}
Suppose $\hat{f}$ is supported in $\abs{\xi} \sim \lambda$. Then 
\begin{equation*}
\norm{(D_t^2 - P_y^2) T_{\lambda}^* U(t) T_{\lambda} f}_{L^2} \lesssim \lambda \norm{f}_{L^2}.
\end{equation*}
\end{lemma}
\begin{proof}
We write $u = T_{\lambda}^* U(t) T_{\lambda} f$ and 
\begin{equation*}
(D_t^2 - P_y^2) u = (D_t + P_y)^2 u - 2 P_y (D_t + P_y) u.
\end{equation*}
Since $u$ and $(D_t + P_y) u$ are supported in $\abs{\xi} \sim \lambda$ and $p$ is of order $1$, we have 
\begin{equation*}
\norm{P_y (D_t + P_y) u}_{L^2} \lesssim \lambda \norm{f}_{L^2}
\end{equation*}
by Lemma \ref{halfwaveapproximate}.

Using \eqref{halfwaveu} and writing $(P_y + L_{x,\xi}) g_{\lambda} = (g_{x,\xi})_{\lambda}$, we get 
\begin{equation*}
(D_t + P_y)^2 u = \iint (U(t) T_{\lambda} f)(x,\xi) (P_y + L_{x,\xi}) (g_{x,\xi})_{\lambda}(y;x,\xi) \,dx\,d\xi.
\end{equation*}
Since 
\begin{align*}
D_x (g_{x,\xi})_{\lambda} &= -D_y (g_{x,\xi})_{\lambda} + (D_x g_{x,\xi})_{\lambda}, \\
D_{\xi} (g_{x,\xi})_{\lambda} &= (y-x) (g_{x,\xi})_{\lambda} + (D_{\xi} g_{x,\xi})_{\lambda},
\end{align*}
the argument in Lemma \ref{wavepacketparameter} gives 
\begin{equation*}
(P_y + L_{x,\xi}) (g_{x,\xi})_{\lambda} = (m_{x,\xi}(z,D_z) g_{x,\xi})_{\lambda} + (L_{x,\xi} m_{x,\xi}(z,D_z) g)_{\lambda}.
\end{equation*}
The Schwartz seminorms of $g_{x,\xi}$ are uniformly bounded in $x,\xi$, and \eqref{mxxiest} shows that the same applies to the Schwartz seminorms of $m_{x,\xi}(z,D_z) g_{x,\xi}$. Lemma \ref{fbisobolevbounded} gives 
\begin{equation*}
\norm{ \iint (U(t) T_{\lambda} f)(x,\xi) (m_{x,\xi}(z,D_z) g_{x,\xi})_{\lambda}(y;x,\xi) \,dx\,d\xi }_{L^2(\mR^n_y)} \lesssim \norm{f}_{L^2}
\end{equation*}
It remains to study the symbol seminorms of $\tilde{m}_{x,\xi}(z,\zeta) = L_{x,\xi} m_{x,\xi}(z,\zeta)$. From \eqref{mxxidef} we obtain 
\begin{multline}
\tilde{m}_{x,\xi}(z,\zeta) = \int_0^1 (1-s) \partial_s^2 \Big[ d_{\xi} p(x,\xi) \cdot D_x p(x+s\lambda^{-1/2} z,\xi+s\lambda^{1/2} \zeta) \\
 - d_x p(x,\xi) \cdot D_{\xi} p(x+s\lambda^{-1/2} z,\xi+s\lambda^{1/2} \zeta) \Big] \,ds. \label{mxxitilde_def}
\end{multline}
Since $\abs{\xi} \sim \lambda$ a computation as in \eqref{mxxi_explicit} shows  
\begin{equation} \label{mxxitildeest}
\abs{\partial_z^{\alpha} \partial_{\zeta}^{\beta} \tilde{m}_{x,\xi}(z,\zeta)} \lesssim \lambda^{1/2} \br{z}^2, \qquad \abs{\zeta} \leq 2,
\end{equation}
and Lemma \ref{fbisobolevbounded} implies 
\begin{equation*}
\norm{ \iint (U(t) T_{\lambda} f)(x,\xi) (\tilde{m}_{x,\xi}(z,D_z) g)_{\lambda}(y;x,\xi) \,dx\,d\xi }_{L^2(\mR^n_y)} \lesssim \lambda^{1/2} \norm{f}_{L^2}.
\end{equation*}
This concludes the proof.
\end{proof}

We move on to the full parametrix, where all the frequencies are added up. The errors resulting from truncation are handled as in \cite[Theorem 4.5]{smith1}.

\begin{lemma} \label{stildeestimates}
$\tilde{S}(t)$ is a bounded operator $H^{\alpha} \to H^{\alpha+1}$, for all $\alpha \in \mR$. If $-1 \leq \alpha \leq 2$ then $(D_t^2 - A(x,D_x)) \tilde{S}(t)$ is a bounded operator $H^{\alpha} \to H^{\alpha}$.
\end{lemma}
\begin{proof}
Since $q_k^{\pm}$ is order $-1$ one has $E_k^{\pm}(t): H^{\alpha} \to H^{\alpha+1}$, and since $E_k^{\pm}(t) g$ is localized near $\abs{\xi} \sim 2^k$ the sum converges in $H^{\alpha+1}$ and $\tilde{S}(t): H^{\alpha} \to H^{\alpha+1}$.

Let now $-1 \leq \alpha \leq 2$. We write 
\begin{multline*}
(D_t^2 - A(x,D_x)) \tilde{S}(t) g = -t \sum_{k < k_0} A(x,D_x) \beta_k(D) g \\
 + \sum_{k \geq k_0} (D_t^2 - A_k(x,D_x))(u_k^+ + u_k^-) + \sum_{k \geq k_0} (A_k(x,D_x) - A(x,D_x))(u_k^+ + u_k^-)
\end{multline*}
and write the last expression as $\tilde{T}_1(t)g + \tilde{T}_2(t)g + \tilde{T}_3(t)g$. It is clear that $\norm{\tilde{T}_1(t) g}_{H^{\alpha}} \lesssim \norm{g}_{H^{\alpha}}$. For $\tilde{T}_2(t)$ we write 
\begin{equation*}
D_t^2 - A_k(x,D_x) = D_t^2 - (P_k^{\pm})^2 + (P_k^{\pm})^2 - A_k(x,D_x).
\end{equation*}
One has 
\begin{equation*}
\norm{\sum_{k \geq k_0} (D_t^2 - (P_k^+)^2) u_k^+}_{H^{\alpha}} \lesssim \norm{g}_{H^{\alpha}}
\end{equation*}
by Lemma \ref{psdofullwaveapproximate}. From \eqref{truncation_est} we obtain 
\begin{equation} \label{pkak_estimate}
(\sqrt{A_k} - p_k^+)\beta_k(\xi) \in S^0_{1,1/2},
\end{equation}
which shows that $(P_k^+)^2 - A_k(x,D_x)$ is of order $1$ and 
\begin{equation*}
\norm{\sum_{k \geq k_0} ((P_k^+)^2 - A_k(x,D_x)) u_k^+}_{H^{\alpha}} \lesssim \norm{g}_{H^{\alpha}}.
\end{equation*}
Similar results hold for $u_k^-$.

For $\tilde{T}_3(t)$ we need to show that 
\begin{equation*}
\Gamma: g \mapsto \sum_{k \geq k_0} (a_k(x) - a(x)) D_{x_i} D_{x_j} u_k^+
\end{equation*}
is bounded on $H^{\alpha}$, where $a_k = a^{ij}_k$ and $a = a^{ij}$. We will write $\beta_j(D) \Gamma = \sum_{k \geq k_0} \Gamma_{jk} D_{x_i x_j} u_k^+$ where 
\begin{equation*}
\Gamma_{jk} = \tilde{\beta}_j(D) (a(x)-a_k(x)) \tilde{\beta}_k(D).
\end{equation*}
Here $\tilde{\beta}_k(\xi)$ are cutoffs to $\abs{\xi} \sim 2^k$ which are $1$ on the frequency support of $u_k^+$. We let $l_0 = l_0(M)$ be an integer so that $\supp(\tilde{\beta}_k(\xi)) \subseteq \{2^{-l_0+k} \leq \abs{\xi} \leq 2^{l_0+k} \}$.

By looking at the supports on the Fourier side, we get 
\begin{equation*}
\Gamma_{jk} = \left\{ \begin{array}{ll} 
\beta_j(D) (a-a_{2k-4l_0}) \tilde{\beta}_k(D), & j \leq k - 3l_0, \\
\beta_j(D) (a-a_{k}) \tilde{\beta}_k(D), & k - 3l_0 \leq j \leq k + 3l_0, \\
\beta_j(D) (a-a_{2j-4l_0}) \tilde{\beta}_k(D), & j \geq k + 3l_0.
\end{array} \right.
\end{equation*}
Using that $\norm{a-a_k}_{L^{\infty}} \lesssim 2^{-k}$ we obtain 
\begin{multline*}
\norm{\beta_j(D) \Gamma g}_{L^2}^2 \lesssim \sum_{k < j-3l_0} 2^{-4j} 2^{2k} \norm{g_k}_{L^2}^2 + \sum_{\abs{j-k} \leq 3l_0} 2^{-2j}2^{2k} \norm{g_k}_{L^2}^2 \\
 + \sum_{k > j+3l_0} 2^{-4k} 2^{2k} \norm{g_k}_{L^2}^2.
\end{multline*}
It follows that $\Gamma$ is bounded on $H^{\alpha}$ with $\abs{\alpha} \leq 1$. We compute 
\begin{equation*}
D_{x_l} \Gamma g = \sum_{k \geq k_0} (D_l a_k) D_{x_i x_j} u_k^+ - D_l a \sum_{k \geq k_0} D_{x_i x_j} u_k^+ + \sum_{k \geq k_0} (a_k-a) D_{x_i x_j x_l} u_k^+
\end{equation*}
and note that this is bounded $H^{\alpha+1} \to H^{\alpha}$ if $\abs{\alpha} \leq 1$, due to the $C^{1,1}$ regularity of $a$ and the argument above. This shows boundedness of $\Gamma$ for $-1 \leq \alpha \leq 2$.
\end{proof}

In the next lemma we correct the value of $\partial_t \tilde{S}(t)$ at $t = 0$.

\begin{lemma}
If $k_0$ is sufficiently large, then 
\begin{equation*}
K = \sum_{k \geq k_0} (i \tilde{R}_k^+(0) + i \tilde{R}_k^-(0) + \frac{1}{2} R_k^+ \beta_k(D) + \frac{1}{2} R_k^- \beta_k(D))
\end{equation*}
has norm $\leq 1/2$ on $H^{\alpha}$. The operator 
\begin{equation*}
\widehat{S}(t) = \tilde{S}(t) (I+K)^{-1}
\end{equation*}
will satisfy $\widehat{S}(t) g|_{t=0} = 0$ and $\partial_t \widehat{S}(t) g|_{t=0} = g$.
\end{lemma}
\begin{proof}
One has $\norm{\tilde{R}_k^{\pm}(0) g}_{L^2} \lesssim 2^{-k} \norm{g}_{L^2}$ by Lemma \ref{halfwaveapproximate}. Also, $R_k^{\pm} \beta_k(D)$ is of order $-1$ by looking at the symbol expansion of $(P_k^{\pm} Q_k^{\pm} - I) \beta_k(D)$ and using \eqref{pkak_estimate}. The norms do not depend on any previous value of $k_0$, and we may choose $k_0$ so large that $K$ will have norm $\leq 1/2$ on $H^{\alpha}$ for $\alpha \in [-\alpha_0,\alpha_0]$, for any $\alpha_0 > 0$. Thus $I+K$ will be invertible on these spaces and the norm of the inverse will be $\leq 2$.

Since $\tilde{S}(0)g = 0$, also $\widehat{S}(0)g = 0$. The derivative is 
\begin{equation*}
\partial_t \tilde{S}(t) g = \sum_{k < k_0} g_k + i \sum_{k \geq k_0} (D_t u_k^+ + D_t u_k^-).
\end{equation*}
We write $D_t u_k^{\pm} = \tilde{R}_k^{\pm}(t) g - P_k^{\pm} u_k^{\pm}$. Since $u_k^{\pm}(0) = \frac{i}{2} Q_k^{\pm} g_k$, we get 
\begin{equation*}
\partial_t \tilde{S}(t) g|_{t=0} = \sum_{k < k_0} g_k + \sum_{k \geq k_0} (i \tilde{R}_k^+(0) g + i \tilde{R}_k^-(0) g + \frac{1}{2} P_k^{+} Q_k^{+} g_k + \frac{1}{2} P_k^{-} Q_k^{-} g_k).
\end{equation*}
It follows that $\partial_t \tilde{S}(t) g|_{t=0} = (I+K)g$.
\end{proof}

It now remains to show that one obtains a full solution operator from $\widehat{S}(t)$ by introducing a correction by solving a Volterra equation. The solution operator for the Volterra equation is the following.

\begin{lemma} \label{lemma:volterra}
Suppose $T(t,s)$ is bounded on $H^{\alpha}$, with norm uniformly bounded by $C$ when $t,s \in [-M,M]$. There is a bounded map $V$ on $L^{\infty}_t H^{\alpha}_x$, with norm bounded by $e^{CM}$, such that for any $F(t,x) \in L^{\infty}_t H^{\alpha}_x$, $G = VF$ solves the equation 
\begin{equation*}
G(t,x) - \int_0^t T(t,s)G(s,x) \,ds = F(t,x).
\end{equation*}
\end{lemma}
\begin{proof}
Define $(VF)(t,x)$ by  
\begin{equation*}
F(t,x) + \sum_{j=1}^{\infty} \int_0^t \!\! \int_0^{s_1} \!\!\!\!\! \cdots \! \int_0^{s_{j-1}} \! T(t,s_1) T(s_1,s_2) \cdots T(s_{j-1},s_j) F(s_j,x) \,ds_j \cdots \,ds_1.
\end{equation*}
It is easy to check that the series converges in $L^{\infty}_t H^{\alpha}_x$ to a solution which satisfies the desired norm estimate.
\end{proof}

As in the outline, we write $\widehat{S}(t,s) = \widehat{S}(t-s)$ as the operator corresponding to $\widehat{S}(t)$ but with initial surface $\{t = s\}$. If $s,t \in [-M,M]$ then $\widehat{S}(t,s)$ has the same properties as $\widehat{S}(t)$, except that $\widehat{S}(t,s) g|_{t=s} = 0$ and $\partial_t \widehat{S}(t,s)g |_{t=s} = g$. We let 
\begin{equation*}
T(t,s) = (D_t^2 - A(x,D_x)) \widehat{S}(t,s), \qquad T(t) = T(t,0),
\end{equation*}
so $T(t,s)$ and $T(t)$ are bounded on $H^{\alpha}$ when $-1 \leq \alpha \leq 2$ by Lemma \ref{stildeestimates}. We let $V$ be the solution operator for the Volterra equation corresponding to $T(t,s)$.

The argument in the outline shows that if 
\begin{equation*}
S(t)g(x) = \widehat{S}(t)g(x) + \int_0^t \widehat{S}(t,s) V(T(s) g(x)) \,ds,
\end{equation*}
then $(D_t^2 - A(x,D_x))S(t)g = 0$ and $S(t)g|_{t=0} = 0$, $\partial_t S(t) g|_{t=0} = g$. One only needs to check that the time derivatives are justified, but this may be done as in \cite{smith1} and will not be needed for stability considerations. This ends the construction of the solution operator.

\section{Stability} \label{sec:stability}

We now proceed to prove the stability part of the result. Let $A(x) = (a^{ij}(x))$ and $B(x) = (b^{ij}(x))$ be two symmetric matrices satisfying \eqref{a_assumptions}, and take $M$ so large that the truncated metrics $(a^{ij}_k)$ and $(b^{ij}_k)$ also satisfy \eqref{a_assumptions}. Also assume that $t$ satisfies \eqref{t_assumptions}.

We write $P_A = P_k^{\pm}$ for the operator at frequency $\lambda = 2^k$ defined in terms of the metric $A$, and similarly $Q_A, R_A, \tilde{R}_A(t)$ etc. The following operators depend norm continuously on the metric.

\begin{lemma} \label{psdostability}
If $\hat{f}$ is supported in $\abs{\xi} \sim \lambda$, then 
\begin{align*}
\norm{(P_A - P_B)f}_{L^2} &\lesssim \lambda \norm{A-B}_{L^{\infty}} \norm{f}_{L^2}, \\
\norm{(Q_A - Q_B)f}_{L^2} &\lesssim \lambda^{-1} \norm{A-B}_{L^{\infty}} \norm{f}_{L^2}, \\
\norm{(R_A - R_B)f}_{L^2} &\lesssim \lambda^{-1} \norm{A-B}_{L^{\infty}} \norm{f}_{L^2}.
\end{align*}
\end{lemma}
\begin{proof}
Consider $h(x,\xi) = (\lambda^{n/2} \chi(\lambda^{1/2} \,\cdot\,) \ast [F(A_k(\,\cdot\,,\xi)) - F(B_k(\,\cdot\,,\xi))]) \tilde{\beta}_{\lambda}(\xi)$, where $\tilde{\beta}_{\lambda}$ is a cutoff to $\abs{\xi} \sim \lambda$ and $F(t) = t^{1/2}$. We wish to show that 
\begin{equation} \label{papb_est}
\abs{\partial_x^{\alpha} \partial_{\xi}^{\beta} h(x,\xi)} \leq C_{M,\alpha,\beta} \norm{A-B}_{L^{\infty}} \lambda^{1-\abs{\beta} + \frac{1}{2} \abs{\alpha}}.
\end{equation}
This will show the estimate for $P_A - P_B$, and the estimate for $Q_A - Q_B$ follows from a similar result with $F(t) = t^{-1/2}$.

In $\partial_x^{\alpha} \partial_{\xi}^{\beta} h$ we let the $x$-derivatives hit the mollifier, which gives the desired growth. We may thus assume that $\alpha = 0$. Each $\xi$-derivative hitting $\tilde{\beta}_{\lambda}(\xi)$ gives $\lambda^{-1}$, so we only need to look at the case when the $\xi$-derivatives hit $F(A_k) - F(B_k)$. We write 
\begin{equation*}
F(A_k) - F(B_k) = \int_0^1 F'(r A_k + (1-r) B_k) (A_k - B_k) \,dr.
\end{equation*}
The matrix $r A_k + (1-r) B_k$ satisfies \eqref{a_assumptions}. Consequently 
\begin{equation*}
\abs{\partial_{\xi}^{\beta} [ (F(A_k) - F(B_k)) \tilde{\beta}_k(\xi) ]} \leq C_{M,\beta} \norm{A-B}_{L^{\infty}} \lambda^{1-\abs{\beta}}
\end{equation*}
as desired.

If $P, Q$ are pseudodifferential operators the symbol of $PQ$ is 
\begin{equation*}
\sigma(PQ) = \sum_{\abs{\alpha} < N} \frac{\partial_{\xi}^{\alpha} p D_x^{\alpha} q}{\alpha!} + \sum_{\abs{\alpha+\beta}=2N} \int e^{-iy \cdot \eta} s_{\alpha \beta} (x,\xi,y,\eta) \,dy\,d\eta
\end{equation*}
where the last terms are oscillatory integrals, and 
\begin{equation*}
s_{\alpha \beta} = \sum_{\gamma \leq \alpha,\beta} c_{\alpha \beta \gamma} \int_0^1 t^{\abs{\alpha+\beta-2\gamma}} (1-t)^{2N-1} \partial_{\xi}^{\alpha+\beta-\gamma} p(x,\xi+t\eta) \partial_x^{\alpha+\beta-\gamma} q(x+ty,\xi) \,dt.
\end{equation*}
Note that $\abs{\alpha+\beta-\gamma} \geq N$. We have 
\begin{equation*}
(R_A - R_B) \beta_{\lambda}(D) = P_A \tilde{\beta}_{\lambda}(D) Q_A \beta_{\lambda}(D) - P_B \tilde{\beta}_{\lambda}(D) Q_B \beta_{\lambda}(D)
\end{equation*}
for a suitable $\tilde{\beta}_{\lambda}$. Suppressing the cutoffs, this has the symbol 
\begin{multline*}
\sigma(R_A - R_B) = p_A q_A - p_B q_B + \sum_{0 < \abs{\alpha} < N} \frac{\partial_{\xi}^{\alpha} p_A D_x^{\alpha} q_A - \partial_{\xi}^{\alpha} p_B D_x^{\alpha} q_B}{\alpha!} \\
 + \sum_{\abs{\alpha+\beta}=2N} \int e^{-iy \cdot \eta} (s^A_{\alpha \beta} (x,\xi,y,\eta) - s^B_{\alpha \beta} (x,\xi,y,\eta)) \,dy\,d\eta
\end{multline*}
The principal symbol is $p_A q_A - p_B q_B = (1-\chi(2^{-k/2} D_x))(1/p_B - 1/p_A)$ since $q = \chi(2^{-k/2} D_x) (1/p)$, and the arguments above show that the corresponding operator has norm $\lesssim \lambda^{-1} \norm{A-B}_{L^{\infty}}$ on $L^2$. It is easy to see that the terms with $0 < \abs{\alpha} < N$ have the same bound. Finally, if $\abs{\gamma} + \abs{\delta} \leq 2n+1$ then 
\begin{equation*}
\abs{\partial_y^{\gamma} \partial_{\eta}^{\delta} \partial_x^{\alpha'} \partial_{\xi}^{\beta'} (s^A_{\alpha \beta} - s^B_{\alpha \beta})(x,\xi,y,\eta)} \leq C \lambda^{-\frac{N}{2}+n+1-\abs{\beta'}-\frac{1}{2}\abs{\alpha'}} \norm{A-B}_{L^{\infty}}.
\end{equation*}
Taking $N$ large enough and using standard estimates for oscillatory integrals gives the $L^2$ bound for $R_A - R_B$.
\end{proof}

Let $g_{x,\xi}^A$ be the Schwartz functions in Lemma \ref{wavepacketparameter}, defined in terms of the metric $A$.

\begin{lemma} \label{wavepacketparameterstability}
The Schwartz norms of $g_{x,\xi}^A - g_{x,\xi}^B$ are $\lesssim \lambda \norm{A-B}_{L^{\infty}}$, uniformly in $x$ and $\xi$.
\end{lemma}
\begin{proof}
It is enough to show that 
\begin{equation} \label{mxxidiff_est}
\abs{\partial_z^{\alpha} \partial_{\zeta}^{\beta} (m_{x,\xi}^A - m_{x,\xi}^B)(z,\zeta)} \lesssim \lambda \norm{A-B}_{L^{\infty}} \br{z}^2, \quad \abs{\zeta} \leq 2.
\end{equation}
This follows from the expression \eqref{mxxidef} for $m_{x,\xi}$, the computation \eqref{mxxi_explicit}, and the estimate \eqref{papb_est}.
\end{proof}

\begin{remark}
It is also true that the Schwartz seminorms of $g_{x,\xi}^A - g_{x,\xi}^B$ are $\lesssim \norm{A-B}_{C^{1,1}}$, which follows by replacing \eqref{papb_est} with the alternate estimate 
\begin{equation*}
\abs{\partial_x^{\alpha} \partial_{\xi}^{\beta} h(x,\xi)} \leq C_{M,\alpha,\beta} \norm{A-B}_{C^{1,1}} \lambda^{1-\abs{\beta} + \frac{1}{2} \max(0,\abs{\alpha}-2)}.
\end{equation*}
However, due to translation along Hamilton flow, there is a loss of one derivative in the stability estimate in any case. Therefore we can afford to lose one derivative in other estimates as well. This results in stability in terms of the $C^{0,1}$ norm of the metric instead of $C^{1,1}$.
\end{remark}

\begin{lemma}
$\norm{(\tilde{R}_A(0) - \tilde{R}_B(0))f}_{L^2} \lesssim \norm{A-B}_{L^{\infty}} \norm{f}_{L^2}$.
\end{lemma}
\begin{proof}
One has 
\begin{multline*}
(\tilde{R}_A(0) - \tilde{R}_B(0))f =\! \frac{i}{2} \Big( \iint (T_{\lambda}(Q_A-Q_B)\beta_{\lambda}(D)f)(x,\xi) (g_{x,\xi}^A)_{\lambda}(y;x,\xi) \,dx\,d\xi \\
 + \iint (T_{\lambda} Q_B \beta_{\lambda}(D)f)(x,\xi) (g_{x,\xi}^A-g_{x,\xi}^B)_{\lambda}(y;x,\xi) \,dx\,d\xi \Big).
\end{multline*}
The result follows from Lemmas \ref{fbisobolevbounded}, \ref{wavepacketparameter}, \ref{psdostability} and \ref{wavepacketparameterstability}.
\end{proof}

From the preceding results, and from the factorization 
\begin{equation*}
(I-K_A)^{-1} - (I-K_B)^{-1} = (I-K_A)^{-1}(K_B-K_A)(I-K_B)^{-1},
\end{equation*}
we see that $(I+K)^{-1}$, the operator which corrects the initial values, depends norm continuously on the metric: 
\begin{equation} \label{boundarycorrectionstability}
\norm{(I+K_A)^{-1} - (I+K_B)^{-1}}_{H^{\alpha} \to H^{\alpha}} \lesssim \norm{A-B}_{L^{\infty}}.
\end{equation}

Next we consider the stability of the flow operator. Here we will lose one derivative to get Lipschitz stability.

\begin{lemma} \label{flowstability}
If $\hat{f}$ vanishes unless $\abs{\xi} \sim \lambda$, then 
\begin{equation*}
\norm{T_{\lambda}^* (U_A - U_B) T_{\lambda} f}_{L^2} \lesssim \lambda \norm{A-B}_{C^{0,1}} \norm{f}_{L^2}.
\end{equation*}
\end{lemma}
\begin{proof}
Fix $t$ and write $\chi_A = \chi_{0,t}^A$, $\chi_B = \chi_{0,t}^B$. One has 
\begin{equation*}
T_{\lambda}^* (U_A - U_B) T_{\lambda} f(y) = \iint [T_{\lambda} f(\chi_A(x,\xi)) - T_{\lambda} f(\chi_B(x,\xi))] g_{\lambda}(y;x,\xi) \,dx\,d\xi.
\end{equation*}
We let $C_r(x,\xi) = r A_{\lambda}(x,\xi) + (1-r) B_{\lambda}(x,\xi)$, and let $\Phi_r = \chi_{C_r} = (x_r,\xi_r)$ be the flow corresponding to the metric $C_r$. Then $(\Phi_r)_{r \in [0,1]}$ is a smooth family of symplectic diffeomorphisms of $T^* \mR^n$, and we have 
\begin{equation*}
T_{\lambda} f(\chi_A(x,\xi)) - T_{\lambda} f(\chi_B(x,\xi)) = \int_0^1 (d_{x,\xi} T_{\lambda} f)(\Phi_r(x,\xi)) \cdot \partial_r \Phi_r(x,\xi) \,dr.
\end{equation*}
Let $h(s) = h(s,r,x,\xi) = (x_r(s), \xi_r(s)/\lambda)$ where $(x_r(0),\xi_r(0)/\lambda) = (x,\xi/\lambda)$. Differentiating the Hamilton equations for $(x_r,\xi_r/\lambda)$ with respect to $r$, and using 
\begin{equation*}
\partial_r p_{C_r} = \chi(\lambda^{-1/2} D_x) \frac{A_{\lambda}(x,\xi)-B_{\lambda}(x,\xi)}{2\sqrt{C_r(x,\xi)}},
\end{equation*}
we obtain 
\begin{equation*}
\abs{(\partial_r h)\edot(s)} \lesssim \abs{\partial_r h(s)} + \norm{A-B}_{C^{0,1}}.
\end{equation*}
Since $\partial_r h(0) = (0,0)$, Gronwall's inequality shows $\abs{\partial_r h(s)} \lesssim \norm{A-B}_{C^{0,1}}$ for $\abs{s} \leq M$. This implies 
\begin{align*}
\abs{\partial_r x_r} &\lesssim \norm{A-B}_{C^{0,1}}, \\
\abs{\partial_r \xi_r} &\lesssim \lambda \norm{A-B}_{C^{0,1}}.
\end{align*}
Using that $\partial_{x_j} T_{\lambda} f = T_{\lambda}(\partial_j f)$ and $\partial_{\xi_j} T_{\lambda} f = \lambda^{-1/2} \tilde{T}_{j,\lambda} f$ where 
$\tilde{T}_{j,\lambda} f(x,\xi) = (f, (iz_j g)_{\lambda}(\,\cdot\,;x,\xi))$, we get 
\begin{multline} \label{uaub_difference}
T_{\lambda}^* (U_A - U_B) T_{\lambda} f(y) = \int_0^1 \iint \Big[ (\partial_r x_r)_j(x,\xi) T_{\lambda} (\partial_j f)(\Phi_r(x,\xi)) g_{\lambda}(y;x,\xi) \\
 + (\partial_r \xi_r)_j(x,\xi) \lambda^{-1/2} \tilde{T}_{j,\lambda} f(\Phi_r(x,\xi)) g_{\lambda}(y;x,\xi) \Big] \,dx\,d\xi\,dr.
\end{multline}
Lemmas \ref{fbisobolevbounded} and \ref{wavepacketparameter} imply the desired estimate.
\end{proof}

\begin{lemma} \label{localizederrorstability}
Suppose $\hat{f}$ is supported in $\abs{\xi} \sim \lambda$. Write $M_A f = (D_t + P_y^A) T_{\lambda}^* U_A T_{\lambda} f$ and $N_A f = (D_t^2 - (P_y^A)^2) T_{\lambda}^* U_A T_{\lambda} f$. Then 
\begin{align}
\norm{(M_A - M_B) f}_{L^2} &\lesssim \lambda \norm{A-B}_{C^{0,1}} \norm{f}_{L^2}, \label{halfwaveerrorstability} \\
\norm{(N_A - N_B) f}_{L^2} &\lesssim \lambda^2 \norm{A-B}_{C^{0,1}} \norm{f}_{L^2}. \label{fullwaveerrorstability}
\end{align}
\end{lemma}
\begin{proof}
From \eqref{halfwaveu} we have 
\begin{multline*}
(M_A - M_B) f = \iint (U_A - U_B) T_{\lambda} f(x,\xi) (g_{x,\xi}^A)(y;x,\xi) \,dx\,d\xi \\
 + \iint U_B T_{\lambda} f(x,\xi) (g_{x,\xi}^A - g_{x,\xi}^B)_{\lambda}(y;x,\xi) \,dx\,d\xi
\end{multline*}
where $g_{x,\xi}^A$ is as in Lemma \ref{wavepacketparameter}. That lemma, the argument in Lemma \ref{flowstability}, and Lemma \ref{wavepacketparameterstability} give \eqref{halfwaveerrorstability}. For the other estimate we write 
\begin{equation*}
D_t^2 - P_y^2 = (D_t + P_y)(D_t + P_y) - 2P_y(D_t + P_y)
\end{equation*}
which gives 
\begin{multline*}
(N_A - N_B) f = [(D_t + P_y^A) M_A - (D_t + P_y^B) M_B] f \\
 - 2 (P_y^A - P_y^B) M_A f - 2 P_y^B (M_A - M_B) f.
\end{multline*}
The last two terms have $L^2$ norms $\lesssim \lambda^2 \norm{A-B}_{C^{0,1}} \norm{f}_{L^2}$ by Lemma \ref{psdostability} and \eqref{halfwaveerrorstability}. The first term has the form 
\begin{multline*}
[(D_t + P_y^A) M_A - (D_t + P_y^B) M_B] f =\!\! \iint (U_A - U_B) T_{\lambda} f(x,\xi) (\tilde{g}_{x,\xi}^A)(y;x,\xi) \,dx\,d\xi \\
 + \iint U_B T_{\lambda} f(x,\xi) (\tilde{g}_{x,\xi}^A - \tilde{g}_{x,\xi}^B)_{\lambda}(y;x,\xi) \,dx\,d\xi
\end{multline*}
where $\tilde{g}_{x,\xi}^A = m_{x,\xi}^A(z,D_z) g_{x,\xi}^A + \tilde{m}_{x,\xi}^A(z,D_z) g$, using the notation in Lemma \ref{psdofullwaveapproximate}. The Schwartz seminorms of $\tilde{g}_{x,\xi}^A$ are $\lesssim \lambda^{1/2}$ by \eqref{mxxitildeest}, and those of $\tilde{g}_{x,\xi}^A - \tilde{g}_{x,\xi}^B$ are $\lesssim \lambda^2 \norm{A-B}_{L^{\infty}}$ by \eqref{mxxidiff_est}, \eqref{mxxitilde_def}, \eqref{papb_est}. Now \eqref{fullwaveerrorstability} follows in the same way as \eqref{halfwaveerrorstability}.
\end{proof}

We proceed to prove stability results for the operators where all the frequencies are summed up.

\begin{lemma} \label{lemma:shatstability}
$\norm{(\widehat{S}_A(t)-\widehat{S}_B(t)) g}_{H^{\alpha}} \lesssim \norm{A-B}_{C^{0,1}} \norm{g}_{H^{\alpha}}$.
\end{lemma}
\begin{proof}
Because of \eqref{boundarycorrectionstability} it is enough to prove the estimate for $\tilde{S}_A(t)-\tilde{S}_B(t)$. But we have 
\begin{equation*}
(\tilde{S}_A(t)-\tilde{S}_B(t))g = \sum_{k \geq k_0} ( (E_{k,A}^+(t) - E_{k,B}^+(t))g + (E_{k,A}^-(t) - E_{k,B}^-(t))g )
\end{equation*}
where 
\begin{multline*}
(E_{k,A}^{\pm}(t) - E_{k,B}^{\pm}(t))g = T_k^* [U_{k,A}^{\pm}(t) - U_{k,B}^{\pm}(t)] T_k(\frac{i}{2} Q_{k,A}^{\pm} \beta_k(D) g) \\
 + T_k^* U_{k,B}^{\pm}(t) T_k(\frac{i}{2} [Q_{k,A}^{\pm} - Q_{k,B}^{\pm}] \beta_k(D) g).
\end{multline*}
By Lemmas \ref{psdostability} and \ref{flowstability} one gets 
\begin{equation} \label{estability}
\norm{(E_{k,A}^{\pm}(t) - E_{k,B}^{\pm}(t))g}_{L^2} \lesssim \norm{A-B}_{C^{0,1}} \norm{g}_{L^2}.
\end{equation}
We sum up these estimates and use frequency localization to obtain 
\begin{equation*}
\norm{(\tilde{S}_A(t)-\tilde{S}_B(t))g}_{H^{\alpha}} \lesssim \norm{A-B}_{C^{0,1}} \norm{g}_{H^{\alpha}}.
\end{equation*}
\end{proof}

The next lemma considers the error $T_A(t) = (D_t^2 - A(x,D_x)) \widehat{S}_A(t)$.

\begin{lemma} \label{lemma:tstability}
$\norm{(T_A(t) - T_B(t)) g}_{H^{\alpha}} \lesssim \norm{A-B}_{C^{0,1}} \norm{g}_{H^{\alpha+1}}$ when $\abs{\alpha} \leq 1$.
\end{lemma}
\begin{proof}
As above, it is enough to consider the operator 
\begin{multline*}
\tilde{T}_A(t)g = (D_t^2 - A(x,D_x)) \tilde{S}_A(t)g = -t \sum_{k < k_0} A(x,D_x) g_k \\
 + \sum_{k \geq k_0} (D_t^2 - A_k(x,D_x))(E_{k,A}^+(t)g + E_{k,A}^-(t)g) \\
 + \sum_{k \geq k_0} (A_k(x,D_x) - A(x,D_x))(E_{k,A}^+(t)g + E_{k,A}^-(t)g).
\end{multline*}
We write the last three terms as $\tilde{T}_{A,j}(t) g$ for $j = 1,2,3$. Clearly 
\begin{equation*}
\norm{(\tilde{T}_{A,1}(t) - \tilde{T}_{B,1}(t)) g}_{H^{\alpha}} \lesssim \norm{A-B}_{C^{0,1}} \norm{g}_{H^{\alpha}}.
\end{equation*}
For the second term we use that 
\begin{equation*}
(D_t^2 - A_k(x,D_x)) E_{k,A}^{\pm}(t) g = (D_t^2 - (P_{k,A}^{\pm})^2) E_{k,A}^{\pm}(t) g + C_{k,A}^{\pm} E_{k,A}^{\pm}(t) g
\end{equation*}
where $C_{k,A}^{\pm} = (P_{k,A}^{\pm})^2 - A_k(x,D_x)$. The first term is just $N_A (\frac{i}{2} Q_{k,A}^{\pm} g_k)$, and 
\begin{equation*}
\norm{N_A(\frac{i}{2}Q_{k,A}^{\pm} g_k) - N_B(\frac{i}{2}Q_{k,B}^{\pm}g_k)}_{L^2} \lesssim 2^k \norm{A-B}_{C^{0,1}} \norm{g}_{L^2}
\end{equation*}
by \eqref{fullwaveerrorstability} and Lemma \ref{psdostability}. Writing 
\begin{equation*}
C_{k,A}^{\pm} - C_{k,B}^{\pm} = (P_{k,A}^{\pm})^2 - (P_{k,B}^{\pm})^2 - (A_k(x,D_x) - B_k(x,D_x)),
\end{equation*}
and using the argument in the end of Lemma \ref{psdostability}, we see that 
\begin{equation*}
\norm{(C_{k,A}^{\pm} - C_{k,B}^{\pm}) \tilde{\beta}_k(D) g}_{L^2} \lesssim 2^{2k} \norm{A-B}_{L^{\infty}} \norm{g}_{L^2}
\end{equation*}
when $\tilde{\beta}_k(\xi)$ is a cutoff to $\abs{\xi} \sim 2^k$. Lemma \ref{localizederrorstability} shows 
\begin{equation*}
\norm{(\tilde{T}_{A,2}(t) - \tilde{T}_{B,2}(t)) g}_{H^{\alpha}} \lesssim \norm{A-B}_{C^{0,1}} \norm{g}_{H^{\alpha+1}}.
\end{equation*}
For the last term we write 
\begin{multline*}
(\tilde{T}_{A,3}(t)-\tilde{T}_{B,3}(t))g = \sum_{\pm} \sum_{k \geq k_0} (A_k(x,D_x) - A(x,D_x))((E_{k,A}^{\pm}(t)-E_{k,B}^{\pm}(t))g) \\
 + \sum_{\pm} \sum_{k \geq k_0} (A_k(x,D_x) - A(x,D_x) - (B_k(x,D_x) - B(x,D_x))) E_{k,B}^{\pm}(t) g
\end{multline*}
The discussion in Lemma \ref{stildeestimates} and \eqref{estability} give the bound for the first term. The second term is handled as in Lemma \ref{stildeestimates}, except that we use 
\begin{equation*}
\norm{a-a_k(x) - (b-b_k(x))}_{L^{\infty}} \leq \norm{A-B}_{C^{0,1}} 2^{-k/2}.
\end{equation*}
One first gets the bound for $-1 \leq \alpha \leq 0$ and then for $\abs{\alpha} \leq 1$ by computing the derivative. The result follows.
\end{proof}

We note that the last estimate holds for the operators $T(t,s)$ with uniform constants when $t,s \in [-M,M]$. The final estimate we need is for the Volterra solution operator.

\begin{lemma} \label{lemma:volterrastability}
$\norm{(V_A - V_B) F}_{L^{\infty}_t H^{\alpha}_x} \lesssim \norm{A-B}_{C^{0,1}} \norm{F}_{L^{\infty}_t H^{\alpha+1}_x}$ when $\abs{\alpha} \leq 1$.
\end{lemma}
\begin{proof}
Recalling the definition of $V$ from Lemma \ref{lemma:volterra}, we need to estimate the $L^{\infty}_t H^{\alpha}_x$ norm of terms of the form 
\begin{multline*}
I(t,x) = \int_0^t \!\! \int_0^{s_1} \!\!\!\!\! \cdots \! \int_0^{s_{j-1}} \! T_B(t,s_1) \cdots T_B(s_{l-1},s_l) (T_A-T_B)(s_l,s_{l+1}) \\ 
 T_A(s_{l+1},s_{l+2}) \cdots T_A(s_{j-1},s_j) F(s_j,x) \,ds_j\cdots\,ds_1
\end{multline*}
Choose $C = C(M)$ such that for $t,s \in [-M,M]$, 
\begin{eqnarray*}
 & \norm{T_A(t,s) g}_{H^{\alpha}} + \norm{T_B(t,s) g}_{H^{\alpha}} \leq C \norm{g}_{H^{\alpha}}, \quad -1 \leq \alpha \leq 2, & \\
 & \norm{(T_A-T_B)(t,s) g}_{H^{\alpha}} \leq C \norm{A-B}_{C^{0,1}} \norm{g}_{H^{\alpha+1}}, \quad -1 \leq \alpha \leq 1. & 
\end{eqnarray*}
Using the $H^{\alpha}$ estimate for each $T_B$, the $H^{\alpha+1} \to H^{\alpha}$ estimate for $T_A-T_B$, and then the $H^{\alpha+1}$ estimate for each $T_A$ gives 
\begin{equation*}
\norm{I(t,\,\cdot\,)}_{H^{\alpha}} \leq \frac{C^j t^j}{j!} \norm{A-B}_{C^{0,1}} \norm{F}_{L^{\infty}_t H^{\alpha+1}_x}.
\end{equation*}
There are $j$ terms of the form $I(t,x)$ at level $j$. It follows that 
\begin{equation*}
\norm{(V_A - V_B) F}_{L^{\infty}_t H^{\alpha}_x} \leq \Big( \sum_{j=1}^{\infty} \frac{j (CM)^j}{j!} \Big) \norm{A-B}_{C^{0,1}} \norm{F}_{L^{\infty}_t H^{\alpha+1}_x}.
\end{equation*}
\end{proof}

For the full solution operator we have 
\begin{multline*}
(S_A - S_B)(t)g = (\widehat{S}_A - \widehat{S}_B)(t)g + \int_0^t (\widehat{S}_A - \widehat{S}_B)(t,s) V_A(T_A(s)g(x)) \,ds \\
 + \int_0^t \widehat{S}_B(t,s) (V_A - V_B) (T_A(s)g(x)) \,ds \\
 + \int_0^t \widehat{S}_B(t,s) V_B( (T_A-T_B)(s)g(x) ) \,ds.
\end{multline*}
Consequently, for $g \in H^{\alpha}$ with $-1 \leq \alpha \leq 1$, 
\begin{equation*}
\norm{(S_A - S_B)(t)g}_{H^{\alpha+1}} \lesssim \norm{A-B}_{C^{0,1}} \norm{g}_{H^{\alpha+1}}
\end{equation*}
by Lemmas \ref{lemma:shatstability}, \ref{lemma:tstability} and \ref{lemma:volterrastability}. This shows \eqref{lipschitz_stability}.

To prove \eqref{uniform_stability}, take $g \in H^{\alpha}$ with $-1 \leq \alpha < 2$, and write $g = g_s + g_r$ with $g_s \in H^2$ and $\norm{g_{r}}_{H^{\alpha}}$ small. Using the triangle inequality and \eqref{lipschitz_stability}, we get 
\begin{align*}
\norm{(S_A(t)-S_B(t))g_s}_{H^3} &\lesssim \norm{g_s}_{H^2}, \\
\norm{(S_A(t)-S_B(t))g_s}_{H^2} &\lesssim \norm{A-B}_{C^{0,1}} \norm{g_s}_{H^2},
\end{align*}
and by interpolation $\norm{(S_A(t)-S_B(t))g_s}_{H^{3-\kappa}} \lesssim \norm{A-B}_{C^{0,1}}^{\kappa} \norm{g_s}_{H^2}$. Thus, given $\varepsilon > 0$, by choosing $\norm{g_r}_{H^{\alpha}}$ and $\delta > 0$ small enough we obtain 
\begin{multline*}
\norm{(S_A(t)-S_B(t))g}_{H^{\alpha+1}} \leq \norm{(S_A(t)-S_B(t))g_s}_{H^{\alpha+1}} \\
 + \norm{S_A(t)g_r}_{H^{\alpha+1}} + \norm{S_B(t)g_r}_{H^{\alpha+1}} < \varepsilon
\end{multline*}
whenever $\norm{A-B}_{C^{0,1}} < \delta$. This ends the proof of Theorem \ref{thm:main}.

\section{Extensions}

Here we give some extensions of Theorem \ref{thm:main}, following Section 4 of \cite{smith1}. The solution operators in each case are obtained from $\widehat{S}(t)$ similarly as in \cite{smith1}, and the stability results follow from the arguments given in Section \ref{sec:stability}. Therefore we will omit the proofs in this section. We emphasize that for each theorem below, the stability proof is constructive.

\begin{thm} \label{thm:main_full}
Assume \eqref{a_assumptions} -- \eqref{t_assumptions}, and let $-1 \leq \alpha \leq 2$. Suppose that $f \in H^{\alpha+1}$, $g \in H^{\alpha}$, and $F \in L^1_t H^{\alpha}_x$. Then there is a unique solution in $C^0_t H^{\alpha+1}_x \cap C^1_t H^{\alpha}_x$ for the problem 
\begin{equation*}
\left\{ \begin{array}{rl}
(D_t^2 - A(x,D_x))u(t,x) &\!\!\!= F(t,x), \\[4pt]
u|_{t=0} &\!\!\!= f, \\[4pt]
\partial_t u|_{t=0} &\!\!\!= g.
\end{array} \right.
\end{equation*}
The solution satisfies 
\begin{equation} \label{fullcauchy_normestimate}
\norm{u(t)}_{H^{\alpha+1}} \lesssim \norm{f}_{H^{\alpha+1}} + \norm{g}_{H^{\alpha}} + \norm{F}_{L^1_t H^{\alpha}_x}.
\end{equation}
Also, let $A = (a^{ij})$, $B = (b^{ij})$ satisfy \eqref{a_assumptions}, and let $u_A$, $u_B$ be the corresponding solutions. If $-1 \leq \alpha < 2$, then for any $\varepsilon > 0$ there is $\delta > 0$ such that 
\begin{equation} \label{uniform_stability_full}
\norm{u_A(t) - u_B(t)}_{H^{\alpha+1}} < \varepsilon \ \text{ whenever }\ \norm{A-B}_{C^{0,1}} < \delta.
\end{equation}
Further, if $-1 \leq \alpha \leq 1$ and $f \in H^{\alpha+2}, g \in H^{\alpha+1}, F \in L^1 H^{\alpha+1}_x$, then 
\begin{equation} \label{lipschitz_stability_full}
\norm{u_A(t) - u_B(t)}_{H^{\alpha+1}} \lesssim \norm{A-B}_{C^{0,1}} (\norm{f}_{H^{\alpha+2}} + \norm{g}_{H^{\alpha+1}} + \norm{F}_{L^1_t H^{\alpha+1}_x}).
\end{equation}
\end{thm}

\begin{remark}
Here we show how the existence and uniqueness part of Theorem \ref{thm:main_full} can be used to prove stability of the map $A \mapsto u(t,\,\cdot\,)$ in Theorem \ref{thm:main}. It is enough to prove \eqref{lipschitz_stability} for $-1 \leq \alpha \leq 1$, and the uniform continuity will follow as in Theorem \ref{thm:main}. Let $A, B$ be two metrics satisfying \eqref{a_assumptions}, and let $g \in H^{\alpha}$. We denote by $u_A$, $u_B$ the solutions to \eqref{main_cauchy} with data $g$, and write $v = u_A - u_B$. Then $v$ satisfies 
\begin{equation*}
\left\{ \begin{array}{rl}
(D_t^2 - A(x,D_x))v &\!\!\!= (A(x,D_x)-B(x,D_x)) u_B, \\[4pt]
v(0) = \partial_t v(0) &\!\!\!= 0,
\end{array} \right.
\end{equation*}
which implies 
\begin{equation*}
\norm{v(t)}_{H^{\alpha+1}} \lesssim \norm{A-B}_{C^{0,1}} (\sup_{\abs{t} \leq M} \norm{u_B(t)}_{H^{\alpha+2}}) \lesssim \norm{A-B}_{C^{0,1}} \norm{g}_{H^{\alpha+1}}.
\end{equation*}
This is the required bound.
\end{remark}

Now we consider the wave equation with divergence form operator 
\begin{equation*}
A^D(x,D_x) u = D_{x_i}(a^{ij}(x) D_{x_j} u).
\end{equation*}

\begin{thm} \label{thm:divergence_full}
Assume \eqref{a_assumptions} -- \eqref{t_assumptions}, and let $-2 \leq \alpha \leq 1$. Suppose that $f \in H^{\alpha+1}$, $g \in H^{\alpha}$, and $F \in L^1_t H^{\alpha}_x$. Then there is a unique solution in $C^0_t H^{\alpha+1}_x \cap C^1_t H^{\alpha}_x$ for the problem 
\begin{equation*}
\left\{ \begin{array}{rl}
(D_t^2 - A^D(x,D_x))u(t,x) &\!\!\!= F(t,x), \\[4pt]
u|_{t=0} &\!\!\!= f, \\[4pt]
\partial_t u|_{t=0} &\!\!\!= g.
\end{array} \right.
\end{equation*}
The solution satisfies \eqref{fullcauchy_normestimate}. Also, solutions satisfy \eqref{uniform_stability_full} and \eqref{lipschitz_stability_full}, if the ranges for $\alpha$ are replaced by $-2 \leq \alpha < 1$ and $-2 \leq \alpha \leq 0$, respectively.
\end{thm}

Finally, we consider the Laplace-Beltrami type operator 
\begin{equation*}
A^L(x,D_x) u = \rho(x)^{-1} D_{x_i}(\rho(x) a^{ij}(x) D_{x_j} u)
\end{equation*}
where $\rho = (\det\,(a^{ij}))^{1/2}$.

\begin{thm} \label{thm:laplacebeltrami_full}
Assume \eqref{a_assumptions} -- \eqref{t_assumptions}, and let $-1 \leq \alpha \leq 1$. Suppose that $f \in H^{\alpha+1}$, $g \in H^{\alpha}$, and $F \in L^1_t H^{\alpha}_x$. Then there is a unique solution in $C^0_t H^{\alpha+1}_x \cap C^1_t H^{\alpha}_x$ for the problem 
\begin{equation*}
\left\{ \begin{array}{rl}
(D_t^2 - A^L(x,D_x))u(t,x) &\!\!\!= F(t,x), \\[4pt]
u|_{t=0} &\!\!\!= f, \\[4pt]
\partial_t u|_{t=0} &\!\!\!= g.
\end{array} \right.
\end{equation*}
The solution satisfies \eqref{fullcauchy_normestimate}. Also, solutions satisfy \eqref{uniform_stability_full} and \eqref{lipschitz_stability_full}, if the ranges for $\alpha$ are replaced by $-1 \leq \alpha < 1$ and $-1 \leq \alpha \leq 0$, respectively.
\end{thm}

%\bibliography{hyperbolic_stability}
%\bibliographystyle{hamsplain}

\providecommand{\bysame}{\leavevmode\hbox to3em{\hrulefill}\thinspace}
\providecommand{\href}[2]{#2}

\end{document}